\documentclass[12pt,a4paper]{amsart}
\usepackage{amsfonts}
\usepackage[latin1]{inputenc}
\usepackage[english]{babel}
\usepackage{amssymb}
\usepackage{amsmath}
\usepackage{amsthm}

\theoremstyle{definition}
\newtheorem{definition}{Definition}

\theoremstyle{plain}
\newtheorem{theorem}[definition]{Theorem}

\newtheorem{lemma}[definition]{Lemma}

\newtheorem*{corollaryA}{Corollary A}
\newtheorem*{theoremA}{Theorem A}

\newcommand{\D}{\mathbb{D}}
\newcommand{\bd}[1]{\partial #1}

\newcommand{\R}{\mathbb{R}}

\newcommand{\dt}{\;\mathrm{d}t}
\newcommand{\ds}{\;\mathrm{d}s}
\newcommand{\dx}{\;\mathrm{d}x}
\newcommand{\dhaus}{\; \mathrm{d}\mathcal{H}}
\newcommand{\dy}{\;\mathrm{d}y}

\newcommand{\haus}{\mathcal{H}}

\newcommand{\Mod}{\mathsf{M}}
\newcommand{\grad}{\nabla}

\newcommand{\norm}[1]{\lVert #1 \rVert}

\newcommand{\cH}{{\mathcal H}}

\DeclareMathOperator*{\osc}{osc}

\newcommand{\h}{\mathcal{H}^{n-1}}

\newdimen\vintkern\vintkern11pt
\def\vint{-\kern-\vintkern\int}

\newcommand{\defeq}{\mathrel{\mathop:}=}

\begin{document}

\title[Exponential integrability of monotone Sobolev functions]{Sharp exponential integrability for traces of monotone
Sobolev functions} 

\author{Pekka Pankka}
\address{Department of Mathematics, University of Michigan, Ann Arbor,
  MI 48104, USA}
\email{pankka@umich.edu}

\author{Pietro Poggi-Corradini}
\address{Department of Mathematics, Cardwell Hall, Kansas State University,
Manhattan, KS 66506, USA}
\email{pietro@math.ksu.edu}

\author{Kai Rajala}
\address{Department of Mathematics and Statistics, P.~O.~Box 35 (MaD), FI-40014 Univ. of Jyv\"askyl\"a, Finland.}
\email{kirajala@maths.jyu.fi}

\thanks{The third author was supported by the Academy of Finland. 
The authors thank the Department of Mathematics at the University of
Michigan where part of this research was conducted.
} 

\begin{abstract} 
We answer a question posed in \cite{Poggi-CorradiniP:Eggype} on
exponential integrability of functions of restricted $n$-energy. We
use geometric methods to obtain a sharp exponential
integrability result for boundary traces of monotone Sobolev functions
defined on the unit ball.  
\end{abstract}

\maketitle


\section{Introduction}

The following result answered a problem of A.~Beurling, mentioned by
J.~Moser in a famous paper \cite{MoserJ:ShafiT}.

\begin{theoremA}[Chang-Marshall (1985), \cite{chang-marshall:1985ajm}]\label{thm:chang-m}
There is a universal constant $C<\infty$ so that if $f$ is analytic in $\D$, $f(0)=0$, and 
\begin{equation}\label{eq:dirichlet}
\int_\D |f^\prime(z)|^2 \,\frac{dA(z)}{\pi} \leq 1, 
\end{equation}
then 
\[
\int_0^{2\pi}\exp\left(|f^\star(e^{i\theta})|^2\right)\,d\theta \leq C,
\]
where $f^\star$ is the trace of $f$ on $\bd\D$, i.e., $f^\star(\zeta)=\lim_{t\uparrow 1}f(t\zeta)$ for $\cH^1$-a.e.~$\zeta\in \bd\D$.
\end{theoremA}

This result is moreover ``sharp'' in the following sense: the Beurling functions, 
\[
B_a(z)\defeq\left(\log\frac{1}{1-az}\right)\left(\log\frac{1}{1-a^2}\right)^{\frac{-1}{2}} \qquad 0<a<1
\]
are analytic in $\D$, satisfy $B_a(0)=0$ and \eqref{eq:dirichlet}, and
have the property that for any given $\alpha>1$, 
one can choose $a$ so that the integral
\[
\int_0^{2\pi}\exp\left(\alpha |B_a(e^{i\theta})|^2\right)\,d\theta 
\]
is as large as desired.

The following is an easy corollary of the Chang-Marshall Theorem.

\begin{corollaryA}\label{cor:harmonic}
There is a universal constant $C<\infty$ so that if $u:\D\rightarrow \R$ is harmonic with $u(0)=0$ and 
\[
\int_\D |\nabla u(z)|^2 \,\frac{dA(z)}{\pi}\leq 1, 
\]
then 
\[
\int_{0}^{2 \pi} \exp\left(u^\star(e^{i \theta})^2\right) \, d\theta \leq C,
\]
where $u^\star$ is the trace of $u$ on $\bd\D$, i.e., 
$u^\star(\zeta)=\lim_{t\uparrow 1}u(t\zeta)$ for $\cH^1$-a.e. $\zeta\in\bd\D$. 
\end{corollaryA}
This can also be shown to be sharp by considering the real parts of the Beurling functions.

In \cite{Poggi-CorradiniP:Eggype} the last two authors generalized the Chang-Marshall theorem 
to quasiregular mappings in all dimensions. They also asked in \cite{Poggi-CorradiniP:Eggype} 
whether Corollary~A also generalizes, perhaps substituting
``harmonic'' with ``$n$-harmonic''. In this note we show that this is
indeed possible. The key concept is that of a monotone Sobolev
function, whose definition we recall below, and which is quite general,
and includes for instance $n$-harmonic functions.

\section{Main results}\label{sec:main}

For a continuous function $u \colon \Omega \to \R$, we define the oscillation of $u$ on a compact set $K \subset \Omega$ by 
\[
\osc_K u = \max_{x,y \in K} |u(x)-u(y)|. 
\]
We say that $u \colon \Omega \to \R$ is \emph{monotone} if $\osc_{\partial B} u = \osc_{\bar B} u$ for all $n$-balls $B$ compactly contained in $\Omega$.

Given $u \colon B^n\to \R$ in the Sobolev space $W^{1,n}(B^n)$, the radial limit
\[
\tilde u(y) = \lim_{r \to 1} u(ry)
\]
exists at $\cH^{n-1}$-a.e. point $y \in S^{n-1}$. We denote by $\tilde u$
the almost everywhere defined trace of $u$. Moreover, we denote
the $L^p$-norm of a $p$-integrable $g:\Omega \to \R^n$ by
$\norm{g}_p=\norm{g}_{\Omega,p}$. The surface measure
$\haus^{n-1}(S^{n-1})$ of the unit sphere $S^{n-1}$ is
$\omega_{n-1}$. The notations $B^n(r)=B^n(0,r)$, $B^n=B^n(1)$ for
$n$-dimensional balls will be used.  


\begin{theorem}
\label{thm:monotone}
There exists a constant $C=C(n)>0$ so that if $u \in W^{1,n}(B^n)$ is a non-constant continuous monotone function such that $u(0)=0$, then  
\begin{equation}
\label{eq:thm_monotone_integrability}
\int_{S^{n-1}} \exp\left(\alpha (|\tilde u(y)|/\norm{\grad u}_n)^{n/(n-1)}\right) \dhaus^{n-1}(y) \le C,
\end{equation}
where 
\begin{equation}
\label{eq:thm_monotone_constant}
\alpha = (n-1) \left(\frac{\omega_{n-1}}{2}\right)^{\frac{1}{n-1}}.
\end{equation}
\end{theorem}


The continuity assumption in Theorem \ref{thm:monotone} is of
technical nature. By a theorem of Manfredi \cite{ManfrediJ:Monot},
so-called weakly monotone functions in $W^{1,n}$ are always continuous
and monotone in the above sense. In general, $W^{1,n}$-functions need
not be continuous.  

Theorem \ref{thm:monotone} is not true without the
monotonicity assumption. Indeed, since the $n$-capacity of a point is
zero, one can construct Sobolev functions $u_i\in W^{1,n}(B^n)$ so
that $u_i(0)=0$, $\norm{\nabla u_i}_n\leq 1$, and $\tilde{u}_i(y)\geq
i$ for every $y \in S^{n-1}$.  


Our method of proof for Theorem \ref{thm:monotone} has a similar
geometric flavor as in 
 \cite{marshall:1989ark} and in
\cite{Poggi-CorradiniP:Eggype}, and the end-game is again to appeal to
 Moser's original one-dimensional proof. However, the so-called ``egg-yolk''
property, which was the hardest part to establish in the two papers
 cited above, can be quickly established in our present case. It might
come as a surprise then that Theorem \ref{thm:monotone} is sharp, as
we will see in Theorem \ref{thm:esimerkki} below, 
as opposed to the situation in \cite{Poggi-CorradiniP:Eggype}. 

A function $u \in W^{1,p}_{\operatorname{loc}}(\Omega)$ is called \emph{$p$-harmonic}, $1<p<\infty$, if 
\[
\int_{\Omega} |\nabla u|^{p-2} \nabla u \cdot \nabla \phi \dx=0 
\]
for every $C^{\infty}$-smooth test function $\phi$ with compact
support in $\Omega$, see \cite{HeinonenJ:nonl}. Since $p$-harmonic
functions satisfy the maximum principle (\cite[6.5]{HeinonenJ:nonl}),
they are, in particular, monotone.  


The next result shows that the constant $\alpha$ in Theorem
\ref{thm:monotone} is sharp.
\begin{theorem}
\label{thm:esimerkki}
Let $\alpha$ be as in Theorem \ref{thm:monotone}. There exists a sequence of $n$-harmonic functions $u_i \in W^{1,n}(B^n)$ satisfying $\norm{\grad u_i}_n \le  1$ and $u_i(0)=0$, so that 
\[
\int_{S^{n-1}} \exp \left(\beta |\tilde u_i(y)|^{n/(n-1)}\right) \dhaus^{n-1}(y) \rightarrow \infty  \quad \text{as }  i \to \infty    
\]
whenever $\beta > \alpha$. 
\end{theorem}



\section{Proof of Theorem \ref{thm:monotone}}
In this section we assume that $u$ satisfies the assumptions of
Theorem \ref{thm:monotone}. Moreover, by considering balls
$B^n(0,1-1/j)$, for $j$ large, and using Fatou's lemma, we may assume that
the function $u$ in Theorem \ref{thm:monotone} is defined in a
neighborhood of the unit ball.  

\begin{lemma}
\label{lemma:monotone_radius_estimate}
There exists a constant $r_{0}=r_{0}(n)>0$ so that if
$M_{0}\defeq\max_{\bar{B}^n(r_{0})}|u|$, then
\[
\int_{\{|u|\leq M_{0}\}} |\nabla u|^n \dx \ge M_{0}^n. 
\]
\end{lemma}
\begin{proof}
For $0<r<1$ let $m\defeq \max_{\bar{B}^n(r)}|u|$ 
and set $v \defeq \min \{|u|,m\}$. By monotonicity, and since $u(0)=0$, $\osc_{S^{n-1}(t)}v =m$ for every 
$t \geq r$. By the Sobolev embedding theorem on spheres, see e.g. \cite[Lemma
  1]{GehringF:Rinqms} or \cite{MostowG:Quacms}, there exists a constant $a$ depending only on
$n$ such that 
\begin{eqnarray*}
\int_{B^n\setminus \bar B^n(r)} |\grad v|^n \dx 
&=& \int_{r}^1 \left(\int_{S^{n-1}(t)}|\grad v|^n \dhaus^{n-1}\right)\dt \\
&\ge& \int_{r}^1 \frac{\left(\osc_{S^{n-1}(t)} v\right)^n}{a t} \dt 
= \frac{m^n}{a} \log \frac{1}{r}. 
\end{eqnarray*}
The claim follows by choosing $r_{0} \defeq \exp(-a)$. 
\end{proof}

Let $\Gamma$ be a family of paths in a domain $\Omega$. The \emph{$n$-modulus} $\Mod_n (\Gamma)$ of $\Gamma$ is defined as follows: 
\[
\Mod_n (\Gamma) = \inf_{\rho} \int_{\Omega} \rho^n \dx, 
\]
where $\rho \colon \Omega \to [0,\infty]$ is \emph{an admissible function} for $\Gamma$, i.e.~a  Borel function satisfying 
\begin{equation}
\label{oik}
\int_{\gamma} \rho \ds \ge 1 
\end{equation}
for every locally rectifiable $\gamma \in \Gamma$. The family of all
paths joining two sets $A, B \subset \bar{\Omega}$ in $\Omega$ is
denoted by $\Delta(A,B;\Omega)$. We say that a given property holds for $n$-almost every path in a path family $\Gamma$ if the property holds for all paths in $\Gamma\setminus \Gamma_0$, where $\Gamma_0$ is a subfamily of $\Gamma$ having $n$-modulus zero.


\begin{lemma}
\label{lemma:monotone_measure_estimate}
For every $r \in (0,1)$, there exists a constant $c=c(n,r)$, so that 
\begin{equation}
\label{aus}
\haus^{n-1}\left(\{ y \in S^{n-1} \colon | u(y)|\ge s\}\right) \le c \exp\left(-\alpha I_M^s(u)\right) 
\end{equation}
for $s \ge M$. Here $\alpha$ is as in \eqref{eq:thm_monotone_constant}, 
$M=M(r,u)=\max_{S^{n-1}(r)} |u|$, and 
\[
I_M^s(u) = \int_M^s \frac{\dt}{\left( \int_{\{ |u| = t\}} |\grad u|^{n-1} \dhaus^{n-1}\right) ^{1/(n-1)}}.
\]

\end{lemma}


\begin{proof}
Fix $r\in (0,1)$ and $s>M=M(r,u)$. Write
\[
E=E_s \defeq  \{ y \in S^{n-1} \colon | u(y)|\ge s\}
\]
and 
\[
U_M \defeq \{ x\in B^n \colon |u(x)|\ge M\}.
\]
Also, here and in what follows we write
\begin{equation}
\label{eq:At}
A_t \defeq  \int_{\{ |u| = t\}} |\grad u|^{n-1} \dhaus^{n-1}.
\end{equation}
The fact that $A_{t}$ is a Borel function of $t$ is standard, see for
instance \cite{evans-gariepy} p.~117.

We construct an admissible function $\rho$ for $\Delta(\partial U_M, E; B^n)$ as follows: Let $I = I_M^s(u)$, and set 
\[
\rho(x) \defeq   \frac{|\grad u(x)|}{I A_{|u(x)|}^{1/(n-1)}} \chi_{U_M}(x). 
\]

Since every path in $\Delta(\partial U_M,E;B^n)$ has a subpath in
$\Delta(\partial U_M,E;U_M)$, it suffices to show that $\rho$ is
admissible for $n$-almost every path in $\Delta(\partial U_M,E; U_M)$,
i.e. that the set of paths where \eqref{oik} fails has $n$-modulus
zero. Recall that, by Fuglede's theorem \cite[Theorem 3]{FugledeB:Extlfc}, 
$u$ is absolutely continuous on $n$-almost every path. So, for $n$-almost every
$\gamma \in \Delta(\partial U_M,E; U_M)$ 
parameterized by arc length $\ell(\gamma)$, we have, by change of variables 
\begin{eqnarray*}
\int_\gamma \rho \ds &=& \int_0^{\ell(\gamma)}  \frac{|\grad
u(\gamma(t))|}{I A_{|u(\gamma(t))|}^{1/(n-1)}} \dt \ge I^{-1}
\int_0^{\ell(\gamma)} \frac{|(u\circ \gamma)'(t)|}{A_{|(u\circ
\gamma)(t)|}^{1/(n-1)}} \dt \\ 
&\ge& I^{-1} \int_M^s \frac{\dt}{A_{t}^{1/(n-1)}} = 1. 
\end{eqnarray*}
Thus $\rho$ is an admissible function for $\Delta(\partial
U_M,E;U_M)$, and so also for $\Delta(B^n(r),E; U_M)$, by the definition of
$n$-modulus. By the coarea formula, cf. \cite{MalyJ:coare}, we have  
\begin{eqnarray*}
\Mod_n(\Delta(B^n(r),E;B^n)) &\le& \int_{U_M} \rho^n \dx = I^{-n} \int_{U_M} \frac{|\grad u(x)|^n}{A_{|u(x)|}^{n/(n-1)}} \dx \\
&=& I^{-n} \int_M^s \int_{\{ |u|=t\}} \frac{|\grad u(y)|^{n-1}}{A_t^{n/(n-1)}} \dhaus^{n-1}(y) \dt \\
&=& I^{-n} \int_M^s \frac{ A_t }{A_t^{n/(n-1)}} \dt = I^{1-n}.
\end{eqnarray*}

By the conformal invariance of  $n$-modulus, taking inversion with
respect to the unit sphere yields  
\begin{eqnarray*}
2 \Mod_n (\Delta(B^n(r),E;B^n)) \ge\Mod_n(\Delta(S^{n-1}(r)\cup
S^{n-1}(1/r),E;\R^n)).  
\end{eqnarray*}
By spherical symmetrization and \cite[Theorem 4]{GehringF:Symrs},
\begin{eqnarray*}
2 I^{1-n} &\ge& \Mod_n(\Delta(S^{n-1}(r)\cup S^{n-1}(1/r),E;\R^n)) \\
&\ge& \omega_{n-1} \left( \log
\frac{c_2}{\haus^{n-1}(E)^{1/(n-1)}}\right)^{1-n}, 
\end{eqnarray*}
where $c_2$ depends only on $n$ and $r$. See
\cite{Poggi-CorradiniP:Eggype} for further details. This implies
\eqref{aus}.  
\end{proof}


\begin{proof}[Proof of Theorem \ref{thm:monotone}]


We will use the following result of Moser \cite[Equation
(6)]{MoserJ:ShafiT}: If  
$\omega:[0,\infty) \to[0,\infty)$ is absolutely continuous and
satisfies $\omega(0)=0$, $\omega' \geq 0$   
almost everywhere, and 
$$
\int_{0}^{\infty} (\omega'(t))^n\dt \leq 1, 
$$
then 
\begin{equation}
\label{A}
\int_0^{\infty} \exp(\omega(t)^{n/(n-1)}-t)\dt \leq C, 
\end{equation}
where $C>0$ depends only on $n$. 
By scaling invariance of \eqref{eq:thm_monotone_integrability}, we may
assume that  
\begin{equation}
\label{E}
\int_{B^n} |\nabla u|^n \dx=1. 
\end{equation}
Moreover, we fix $r=r_{0}$ and $M=M_{0}$ as in Lemma
\ref{lemma:monotone_radius_estimate}. Then, in particular, $M<1$.  

By the Cavalieri principle, 
\begin{eqnarray*}
&& \int_{S^{n-1}} \exp\left(\alpha |u(x)|^{n/(n-1)}\right) \dhaus^{n-1}(x) \\
&& \quad = \omega_{n-1} + \frac{\alpha n}{n-1} \int_0^\infty
s^{1/(n-1)} \haus^{n-1}(E_s) \exp\left(\alpha s^{n/(n-1)}\right) \ds,  
\end{eqnarray*}
where 
\[
E_s =\{y \in S^{n-1}: |u(y)|\geq s \}. 
\]
Then, by Lemma \ref{lemma:monotone_measure_estimate}, it suffices to bound 
\begin{equation}
\label{star}
\int_0^{||u||_{\infty}} s^{1/(n-1)}\exp\left(\alpha(s^{n/(n-1)}-I_M^s(u))\right) \ds, 
\end{equation}
where $||u||_{\infty}=\max_{y\in S^{n-1}}|u(y)|$, and $I_M^s(u)=0$ for $0<s<M$. We define a function $\psi:[0,\infty) \to [0,\infty)$, 
\begin{eqnarray*}
&&\psi(s)=\left\{ \begin{array}{ll} \mu s, & 0<s<M\\ 
\alpha I_M^s(u) + \mu M, & M \leq s \leq ||u||_{\infty}\\ 
\alpha I_M^{||u||_{\infty}}(u) + \mu M, & s > ||u||_{\infty},  
\end{array}
\right.
\end{eqnarray*}
where 
\begin{equation}
\label{D}
\mu = \alpha \Big( \frac{M}{\int_{\{|u|\leq M\}}|\nabla u|^{n}\dx} \Big)^{1/(n-1)}. 
\end{equation}
Then, by Lemma \ref{lemma:monotone_radius_estimate}, $\mu M \leq \alpha$, and thus we may consider  
\begin{equation}
\label{twostar}
\int_0^{||u||_{\infty}} s^{1/(n-1)} \exp (\alpha s^{n/(n-1)}-\psi(s))\ds   
\end{equation}
instead of \eqref{star}. We define $\phi$ by $\phi(y)=\psi^{-1}(y)$ for $0<y<||\psi||_{\infty}$, and $\phi(y)=||u||_{\infty}$ for $y \geq ||\psi||_{\infty}$. Then, changing variables $y=\psi(s)$ in \eqref{twostar} yields 
\begin{equation}
\label{tristar}
\int_0^{\infty} \exp(\alpha \phi(y)^{n/(n-1)}-y)\phi'(y)\phi(y)^{1/(n-1)}\dy. 
\end{equation}
Integrating by parts, we then have that \eqref{tristar} equals $C_1(n)+C_2(n)T$,\[
T= \int_0^{\infty} \exp( (\alpha^{(n-1)/n}\phi(y))^{n/(n-1)}-y  ) \dy. 
\]
Now, since $\phi$ is absolutely continuous and increasing, and $\phi(0)=0$, Theorem \ref{thm:monotone} follows from Moser's result \eqref{A} if we can show that 
\begin{equation}
\label{F}
\int_0^{\infty} (\alpha^{(n-1)/n}\phi'(y))^n\,dy \leq 1. 
\end{equation}

We have 
\begin{eqnarray*}
\alpha^{(n-1)/n}\phi'(y)=\left\{ \begin{array}{ll} 
\alpha^{(n-1)/n}\mu^{-1}, & 0<y<\mu M  \\ 
\alpha^{-1/n} A_{\phi(y)}^{1/(n-1)}, & \mu M < y < \norm{\psi}_\infty \\
0, & y > \norm{\psi}_\infty,  
\end{array}
\right.
\end{eqnarray*}
where $A_{\phi(y)}$ as in \eqref{eq:At}. Hence, 
\begin{equation}
\label{B}
\alpha^{n-1}\int_0^{\infty} \phi'(y)^n \dy = \alpha^{n-1}\mu^{1-n}M + \alpha^{-1}\int_{\mu M}^{||\psi||_{\infty}}A_{\phi(y)}^{n/(n-1)} \dy. 
\end{equation}
By our choice of $\mu$, the first term equals $\int_{\{ |u|\leq M\}} |\nabla u|^n \dx$. Also, by changing variables $\phi(y)=s$ in the right hand integral, and using the coarea formula, we have 
\begin{equation}
\label{C}
\begin{split}
\alpha^{-1} \int_{\mu M}^{\norm{\psi}_{\infty}} A_{\phi(y)}^{n/(n-1)} \dy &= \int_{\mu M}^{\norm{\psi}_\infty} A_{\phi(y)} \phi'(y) \dy\\
&= \int_M^{\norm{u}_\infty} A_s \ds = \int_{\{ |u|\ge M\} } |\grad u|^n \dx.
\end{split}
\end{equation}
Combining \eqref{B}, \eqref{C}, \eqref{D} and \eqref{E} then yields \eqref{F}. The proof is complete. 

\end{proof}


\section{Proof of Theorem \ref{thm:esimerkki}}

Fix $\beta > \alpha$. For notational convenience, we consider first functions in $B^n(e_n,1)$ instead of $B^n$. Fix $2 \leq i \in \mathbb{N}$, and denote $\varepsilon=\varepsilon_i=i^{-1}$. 
Define $v=v_i \colon B^n(-\varepsilon e_n,2+\varepsilon) \to \mathbb{R}$, 
\[
v(x)= - \log|x+\varepsilon e_n|. 
\]
Then $v$ is $n$-harmonic in $B^n(-\varepsilon e_n, 2+\varepsilon)\setminus \{-\varepsilon e_n\}$. We first show that
\begin{equation}
\label{pregrad}
\int_{B^n(e_n,1)}|\nabla v|^n \dx \leq \frac{\omega_{n-1}}{2} \log \frac{2+\varepsilon}{\varepsilon}.
\end{equation}

Clearly, 
\[
\int_{B^n(e_n,1)} |\nabla v |^n \dx \leq \frac{1}{2} \int_A |\nabla v|^n \dx, 
\]
where 
\[
A = B^n(-\varepsilon e_n,2+\varepsilon) \setminus \bar B^n(-\varepsilon e_n,\varepsilon).
\]
Since 
\[
|\nabla v(x)|^n= |x+\varepsilon e_n|^{-n},
\]
we have 
\[
\frac{1}{2} \int_A |\nabla v|^n \dx = \frac{1}{2}\int_{B^n(0,2+\varepsilon)\setminus \bar B^n(0,\varepsilon)} |x|^{-n}\,dx= \frac{\omega_{n-1}}{2} \log \frac{2+\varepsilon}{\varepsilon}. 
\]
Hence \eqref{pregrad} holds. 

To study exponential integrability of $v$, set
\begin{equation}
\label{varabeta}
\gamma = \beta \left( \frac{\omega_{n-1}}{2} \log \frac{2+\varepsilon}{\varepsilon} \right)^{1/(1-n)}
\end{equation}
and $\tau = \gamma/(n-1)$.

By the choice of $\gamma$, \eqref{pregrad}, and the Cavalieri principle, 
\begin{equation}
\label{cava}
\begin{split}
&\int_{S^{n-1}(e_n,1)}\exp(\beta \left( |v|/\norm{\grad v}_n \right)^{n/(n-1)}) \dhaus^{n-1} \\
&\quad \ge \omega_{n-1} +  \frac{n\gamma}{n-1}\int_{0}^{\infty} \h (E_s) s^{1/(n-1)} \exp(\gamma s^{n/(n-1)})\,ds, 
\end{split}
\end{equation}
where 
\[
E_s = \{ x \in S^{n-1}(e_n,1) \colon |v(x)|\geq s\}. 
\] 
Since  
\begin{eqnarray*}
E_s &=& S^{n-1}(e_n,1) \cap \bar B^n(-\varepsilon e_n, \exp(-s))\\
&& \cup \ S^{n-1}(e_n,1)\setminus B^n(-\varepsilon e_n, \exp(s)),
\end{eqnarray*}
we have
\begin{equation}
\label{lev}
\h (E_s)  \ge C(n) (\exp(-s))^{n-1} = C(n)\exp((1-n)s)
\end{equation}
for $0 \leq s \leq  \log \left(1/(2\varepsilon)\right)$.

Combining \eqref{cava} and \eqref{lev} yields 
\begin{eqnarray*}
&& \frac{1}{C(n)}\int_{S^{n-1}(e_n,1)} \exp(\gamma |v|^{n/(n-1)}) \dhaus^{n-1}\\
&&\quad \geq \frac{n\gamma}{n-1} \int_0^{\log(1/(2\varepsilon))} s^{1/(n-1)}\exp\left(\gamma s^{n/(n-1)}+(1-n)s\right)\ds \\ 
&&\quad = n\tau \int_0^{\log(1/(2\varepsilon))} s^{1/(n-1)}\exp\left((n-1) (\tau s^{n/(n-1)}-s)\right)\ds \\
&&\quad = \int_0^{\log(1/(2\varepsilon))} (n\tau s^{1/(n-1)}-(n-1)) \exp\left((n-1) (\tau s^{n/(n-1)}-s)\right)\ds \\
&& \qquad + (n-1) \int_0^{\log(1/(2\varepsilon))} \exp\left((n-1) (\tau s^{n/(n-1)}-s)\right)\ds \\
&& 
\label{kraa}
\quad \ge \exp\left((n-1) \left(\tau \left(\log \frac{1}{2\varepsilon}\right)^{n/(n-1)}-\log\frac{1}{2\varepsilon}\right)\right) - 1. 
\end{eqnarray*}

Since
\[
\left( \log \frac{2+\varepsilon}{\varepsilon} \right)^{1/(1-n)} \left(\log \frac{1}{2\varepsilon}\right)^{1/(n-1)} \ge 1-\delta(\varepsilon), 
\]
where $\delta(\varepsilon)\to 0$ as $\varepsilon \to 0$, we have
\begin{equation}
\label{preexponential}
\begin{split}
& \int_{S^{n-1}(e_n,1)} \exp(\gamma |v|^{n/(n-1)}) \dhaus^{n-1}
\quad \ge C(n) \varepsilon^{-T} -C(n), 
\end{split}
\end{equation}
where 
\[
T = (\beta -\alpha)(2/\omega_{n-1})^{1/(n-1)}-\delta'(\varepsilon), 
\]
and $\delta'(\varepsilon)\to 0$ when $\varepsilon \to 0$. 

To prove Theorem \ref{thm:esimerkki}, we consider the sequence $u_i \colon \bar B^n \to \R$,
\[
u_i(x)=v_i(x+e_n) - v_i(e_n),
\]
where $v_i(e_n) = -\log(1-\varepsilon_i) \le \log 2$ for all $i$. We fix $M$ such that 
\[
\beta' = \beta \left(\frac{M-\log 2}{M}\right)^{n/(n-1)} > \alpha.
\]
Set also $E_i=\{ y \in S^{n-1}(e_n,1) \colon |v_i(y)|\ge M\}$. Then
\[
\beta |v_i(y)-v_i(e_n)|^{n/(n-1)} \ge \beta' |v_i(y)|^{n/(n-1)}
\]
on $E_i$ for every $i$.
Thus 
\[
\begin{split}
& \int_{S^{n-1}} \exp(\beta (|u_i|/\norm{\grad u_i}_n)^{n/(n-1)}) \dhaus^{n-1} \\
& \quad = \int_{S^{n-1}(e_n,1)} \exp(\beta (|v_i(y)-v_i(e_n)|/\norm{\grad v_i}_n)^{n/(n-1)}) \dhaus^{n-1}(y)\\
& \quad \ge \int_{E_i} \exp(\beta' (|v_i(y)|/\norm{\grad v_i}_n)^{n/(n-1)}) \dhaus^{n-1}(y)\\
& \quad \ge \int_{S^{n-1}(e_n,1)} \exp(\beta' (|v_i(y)|/\norm{\grad v_i}_n)^{n/(n-1)}) \dhaus^{n-1}(y) \\
& \qquad - \omega_{n-1} \exp(\beta' (M/\norm{\grad v_i}_n)^{n/(n-1)}).
\end{split}
\]
Since $\beta'>\alpha$ and $\varepsilon_i=i^{-1}$ in \eqref{preexponential}, the claim now follows from \eqref{preexponential}.

\bibliographystyle{abbrv}

\end{document}